\documentclass[reqno,12pt]{amsart}

\usepackage{amsmath}
\usepackage{amssymb}
\usepackage{amsthm}
\usepackage[english]{babel}
\newtheorem{theorem}{Theorem}

\newtheorem{proposition}{Proposition}

\newtheorem{corollary}{Corollary}

\begin{document}

\title[Estimates for the concentration functions]
{Estimates for the concentration functions of weighted sums of independent random variables}

\author[Yu.S. Eliseeva]{Yu.S. Eliseeva$^{1}$} 
\author[A.Yu. Zaitsev]{A.Yu. Zaitsev$^{1,2}$}

\email{pochta106@yandex.ru}
\address{St.~Petersburg State University\bigskip}

\email{zaitsev@pdmi.ras.ru}
\address{St.~Petersburg Department of Steklov Mathematical Institute
\\
Fontanka 27, St.~Petersburg 191023, Russia\\
and St.~Petersburg State University}

\begin{abstract}{Let $X,X_1,\ldots,X_n$ be independent identically distributed random variables  The paper deals with the question about the behavior of the concentration function 
of the random variable $\sum_{k=1}^{n}a_k X_k$ according to 
the arithmetic structure of coefficients~$a_k$.
Recently the interest to this question has increased significantly due to the study of distributions of eigenvalues of random matrices.
In this paper we formulate 
and prove some refinements of the results of Friedland and Sodin (2007) 
and~Rudelson and Vershynin (2009).}
\end{abstract}

\keywords {concentration functions, inequalities,
the Littlewood--Offord problem, sums of independent random variables}

\subjclass {Primary 60F05; secondary 60E15, 60G50}

\maketitle

\footnotetext[1]{Research supported by grant RFBR 10-01-00242.}
\footnotetext[2]{Research supported by grant RFBR 11-01-12104  and by the Program of
Fundamental Researches of Russian Academy of Sciences ``Modern
Problems of Fundamental Mathematics''. }

\section{Introduction}

Let $X,X_1,\ldots,X_n$ be independent identically distributed (i.i.d.) random variables 
with the common distribution $F=\mathcal L(X)$.
The L\'evy concentration function of a random variable $X$ is defined by the equality
$$Q(F,\lambda)=\sup_{x\in\mathbf{R}}F\{[x,x+\lambda]\}, \quad \lambda>0.$$
 Let $a=(a_1,\ldots,a_n)\in \mathbf{R}^n$.
This paper deals with the question about the behavior of the concentration function of the random variable $\sum_{k=1}^{n}a_k X_k$ according to the arithmetic structure of coefficients~$a_k$.
Recently the interest to this question has increased significantly due to the study of distributions of eigenvalues of random matrices (see, for instance, Nguyen and Vu \cite{Nguyen and Vu}, Rudelson and Vershynin \cite{Rudelson and Vershynin08},
\cite{Rudelson and Vershynin}, Tao and Vu \cite{Tao and Vu}, \cite{Tao and Vu2}).
The authors of papers mentioned above have called this question the Littlewood--Offord problem.

In the sequel, $F_a$ is the distribution of the sum
$\sum_{k=1}^{n}a_k X_k$, and $G$ is the distribution of the symmetrized random variable
$\widetilde{X}=X_1-X_2$.
Let \begin{equation} \label{0}M(\tau)=\tau^{-2}\int_{|x|\leq\tau}x^2
\,G\{dx\}+\int_{|x|>\tau}G\{dx\}=\mathbf{E}
\min\big\{{\widetilde{X}^2}/{\tau^2},1\big\}, \quad \tau>0.\end{equation}
 The symbol $c$ will be used for absolute positive constants. Note that $c$ can be different in different (or even in the same) formulas.
We will write $A\ll B$ if $A\leq c B$. Also we will write $A\asymp
B$ if $A\ll B$ and $B\ll A$. For~${x=(x_1,\dots,x_n )\in\mathbf R^n}$
we will denote $\|x\|^2= x_1^2+\dots +x_n^2$
and
$\|x\|_\infty= \max_j|x_j|$.

The elementary properties of concentration functions are well studied (see, for instance,
 \cite{Arak and Zaitsev}, \cite{Hengartner and Theodorescu},
\cite{Petrov}).
In particular, it is obvious that $Q(F,\mu)\le (1+\lceil
\mu/\lambda\rceil)\,Q(F,\lambda)$, for any $\mu,\lambda>0$, where $\lceil
x\rceil$ is the integer part of a number~$x$.
Hence, \begin{equation}\label{8a}
Q(F,c\lambda)\asymp\,Q(F,\lambda)
\end{equation}
and\begin{equation}\label{8j}
 \hbox{if }Q(F,\lambda)\ll K,\hbox{ then }Q(F,\mu)\ll K(1+\mu/\lambda).
\end{equation}

Moreover, for any distribution $F$, the classical Ess\'een inequalities hold (\cite{Esseen66}, see also
\cite{Hengartner and Theodorescu} and \cite{Petrov}):
\begin{equation} \label{1}
\lambda\int_{0}^{\lambda^{-1}} {|\widehat{F}(t)|^2
\,dt} \ll Q(F,\lambda) \ll
\lambda\int_{0}^{\lambda^{-1}} {|\widehat{F}(t)|
\,dt},\quad \lambda>0,
\end{equation}
where $\widehat{F}(t)$ is the characteristic function of the corresponding random variable.
Note that upper and lower bounds in \eqref{1} may have different orders. This is due to the presence of the second power of $|\widehat{F}(t)|$ in the left-hand side of \eqref{1}. In the general case both inequalities \eqref{1}
have optimal orders. However, if we assume additionally that the distribution $F$ is symmetric and its characterictic function is non-negative for all~$t\in\mathbf R$, then we have the lower bound:
\begin{equation} \label{1a}Q(F,\lambda)\gg
\lambda\int_{0}^{\lambda^{-1}} {\widehat{F}(t) \,dt}
\end{equation}
and, therefore,
\begin{equation} \label{1b}
Q(F,\lambda)\asymp\lambda\int\limits_{0}^{\lambda^{-1}}
{\widehat{F}(t) \,dt}
\end{equation} (see \cite{Arak and Zaitsev}, Lemma~1.5 of Chapter II).
 The use of relation \eqref{1b} will allow us to simplify the arguments of Friedland and Sodin
 \cite{Friedland and Sodin}  and Rudelson and Vershynin~\cite{Rudelson and Vershynin} which were applied to consider the Littlewood--Offord problem.

We recall now the well-known Kolmogorov--Rogozin inequality \cite{Rogozin}
(see \cite{Arak and
Zaitsev}, \cite{Hengartner and Theodorescu}, \cite{Petrov}).
\medskip

\begin{proposition}\label{thKR} Let $Y_1,\ldots,Y_n$ be independent random variables with the distributions $W_k=\mathcal{L}(Y_k)$.
Let $\lambda_1,\ldots,\lambda_n$ 
be positive numbers, $\lambda_k \leq \lambda$ $(k=1,\ldots,n)$. Then
\begin{equation} \label{2}
Q\Big(\mathcal{L}\Big(\sum_{k=1}^{n}Y_k\Big),\lambda\Big)\ll\lambda\,
\Big(\sum_{k=1}^{n}\lambda_k^2\,\big(1-Q(W_k,\lambda_k)\big)\Big)^{-1/2}.
\end{equation}
\end{proposition}

Ess\'een \cite{Esseen} (see \cite{Petrov}, Theorem 3 of Chapter III)
improved this result. He has shown that the following statement is true.

\begin{proposition}\label{thE}Under the conditions of Proposition $\ref{thKR}$ we have
\begin{equation} \label{3}
Q\Big(\mathcal{L}\Big(\sum_{k=1}^{n}Y_k\Big),\lambda\Big)\ll\lambda
\,\Big(\sum_{k=1}^{n}\lambda_k^2\,
M_k(\lambda_k)\Big)^{-1/2},
\end{equation}
 where $M_k(\tau)=\mathbf{E}\,
\min\big\{{\widetilde{Y_k}^2}/{\tau^2},1\big\}$.\end{proposition}

Also one can find the improvements of \eqref{2} and \eqref{3} in \cite{Arak}, \cite{Arak and Zaitsev}, \cite{Bretagnolle},
\cite{Kesten}, \cite{Miroshnikov and Rogozin} and  \cite{Nagaev and Hodzhabagyan}.
 \medskip

The problem of estimating the concentration function of weighted sums
$\sum_{k=1}^{n}a_k X_k$
under different conditions on the vector $a\in \mathbf{R}^n$ and distributions of summands
was considered in \cite{Friedland and Sodin},
\cite{Nguyen and Vu},  \cite{Rudelson and Vershynin08}, \cite{Rudelson and
Vershynin}, \cite{Tao and Vu} and \cite{Tao and Vu2}. In this paper we formulate 
and prove some refinements of the results \cite{Friedland and Sodin}
and~\cite{Rudelson and Vershynin}.

In order to be able to compare the results of \cite{Friedland and Sodin} and \cite{Rudelson and Vershynin}, we formulate them
using the common notation.

Friedland and Sodin \cite{Friedland and Sodin} have simplified the arguments of Rudelson and Vershynin \cite{Rudelson and Vershynin08} and obtained the following result.

\begin{proposition}\label{thFS} {Let $X, X_1,\ldots,X_n$ be i.i.d. random variables such that
$Q(\mathcal{L}(X),2)\leq1-p$, where $p>0$, and let $a=(a_1,\ldots,a_n)\in
\mathbf{R}^n$. If, for some $D\ge\cfrac{1}{2\,\|a\|_{\infty}}$ and
$\alpha>0$,
\begin{equation} \label{3b}
\|\,t a-m\|\geq\alpha,\ \hbox{ for all \ $m\in \mathbf Z^n$  and \
$t\in\Big[\cfrac{1}{2\,\|a\|_\infty}\,,D\Big]$},
\end{equation}
 then
\begin{equation} \label{3a}
Q(F_a, 1/D)\ll\cfrac{1}{\|a\|D\sqrt p }+\exp\big(-c\, p^2 \alpha^2\big).
\end{equation}}\end{proposition}

In Theorem 3.3 of \cite{Friedland and Sodin}, the statement of Proposition
\ref{thFS}
was formulated and proved in a weakened form. There was~$p^2$ instead of $p$ in the right-hand side of  inequality~\eqref{3a}. However, the possibility to replace $p^2$ by ${p}$ in the result of~\cite{Friedland and Sodin}  follows easily from elementary properties of the concentration function. 
This was observed, for example, in~\cite{Rudelson and Vershynin} (see Proposition \ref{thRV}).

Furthermore, in \cite{Friedland and Sodin}, it was assumed that  $0<D<1$. Moreover, there stands $Q(F_a,1)$ instead of $Q(F_a,1/D)$ in the left-hand side of inequality \eqref{3a}. However, the quantity $Q(F_a,1)$ is, generally speaking, essentially less than
$Q(F_a,1/D)$ for ${0<D<1}$, since then ${1}/{D}>1$.
Nevertheless, the result of Friedland and Sodin \cite{Friedland and Sodin} with
$D=1$ implies inequality \eqref{3a} for any $D>0$ and with $Q(F_a,
1/D)$ instead of $Q(F_a,1)$ as simple as 
Corollary~\ref{c1} is derived below from Theorem~\ref{th1}.

Note that for $|\,t|\leq \cfrac{1}{2\,\|a\|_{\infty}}$ we have
\begin{equation}\label{4s}
\big(\hbox{dist}(ta,\mathbf{Z}^n)\big)^2=\sum_{k=1}^{n} \min_{m_k \in
\mathbf{Z}} |\,ta_k - m_k|^2 =
\sum_{k=1}^{n}|\,ta_k|^2 = t^2\,\|a\|^2, \end{equation} where
$$\hbox{dist}(ta,\mathbf{Z}^n)= \min_{m \in \mathbf{Z}^n}\|\,ta - m\|.
 $$
Therefore, the assumption $|\,t|\geq\cfrac{1}{2\,\|a\|_{\infty}}$
under condition~\eqref{3b} is natural. If $D=\cfrac{1}{2\,\|a\|_{\infty}}$, then condition~\eqref{3b} holds formally for $\alpha=\cfrac{\|a\|}{4\,\|a\|_{\infty}}$. This follows from \eqref{4s}.
Moreover, for $D\ge\cfrac{1}{2\,\|a\|_{\infty}}$,
the quantity~$\alpha$ involved in condition~\eqref{3b} can not be more than
$\cfrac{\|a\|}{4\,\|a\|_{\infty}}$.
\bigskip

The one-dimensional version of multidimensional Theorem 3.3 of Rudelson and Vershynin
\cite{Rudelson and Vershynin} is formulated as follows.

\begin{proposition}\label{thRV} {Let $X, X_1,\ldots,X_n$ be i.i.d. random variables  
such that $Q(\mathcal{L}(X),2)\leq 1-p$, where $p>0$.
Let $\alpha, D>0$  and  $\gamma\in(0,1)$. Assume that \begin{equation}
\label{4d}\|\,ta-m\| \geq \min\{\gamma t\,\|a\|,\alpha\}, \quad\hbox{for all}\  m\in\mathbf{Z}^n \ \hbox{and }\  t\in[0,D].\end{equation}
Then
\begin{equation} \label{4a}Q\Big(F_a,\cfrac{1}{D}\Big) \ll
\cfrac{1}{\gamma D\,\|a\|\,\sqrt{p}}+\exp(-2\,p\,\alpha^2).\end{equation} }
\end{proposition}

The statement of Proposition~\ref{thRV} is formulated in \cite{Rudelson and Vershynin}
in a weakened form. In Theorem~3.3 \cite{Rudelson and Vershynin} it
is assumed that $\|a\|\ge1$, and the factor~$\|a\|$ is absent in the denominator of the fraction in the right-hand side of  inequality~\eqref{4a}.
However, the result of Proposition~\ref{thRV} 
follows easily from the statement of Proposition~\ref{thRV} with $\|a\|=1$.
One should just apply this statement to the vector $b=a/\|a\|$.
Moreover, there is the unnecessary assumption $\mathbf{E} \,X=0$  in the formulation of 
Theorem~3.3 of~\cite{Rudelson and Vershynin}.

 It is clear that if
 \begin{equation} \label{4b}
0<D\le D(a)=D_{\alpha,\gamma}(a)=
\inf\big\{t>0:\hbox{dist}(ta,\mathbf{Z}^n)\leq
\min\{\gamma\|\,ta\|,\alpha\}\big\},
\end{equation}then  condition \eqref{4d} holds.
Rudelson and Vershynin \cite{Rudelson and Vershynin} called the quantity~$D(a)$
the essential least common denominator of the vector ${a\in\mathbf{R}^n}$.

Finally, we have to mention that the real formulation of Theorem 3.3 of~\cite{Rudelson and Vershynin} is in fact a consequence of inequality \eqref{4a} which follows from \eqref{4a} by using  relations
\eqref{8j} and~\eqref{4b}.
 \medskip

Now we formulate the first main result of this paper.

\begin{theorem}\label{th1} {Let $X,X_1,\ldots,X_n$ be i.i.d. random variables. Let
$a=(a_1,\ldots,a_n)\in \mathbf{R}^n$
and, for some $\alpha>0$,  condition \eqref{3b} holds with $D=1$, i.e. \begin{equation} \label{5b}
\|\,t a-m\|\geq\alpha\ \hbox{ for all $m\in \mathbf Z^n$  and \
$t\in\Big[\cfrac{1}{2\,\|a\|_\infty},1\Big]$}.
\end{equation}
 Then
$$ Q(F_a,1)\ll \frac{1}{\|a\|\sqrt{M(1)}}+\exp\big(-c\,\alpha^2 M(1)\big),$$
where the quantity $M(1)$ is defined by formula \eqref{0}.}\end{theorem}
\medskip

Now we formulate what follows from Theorem\/ $\ref{th1}$ under the conditions of Proposition~\ref{thFS}.

\begin{corollary}\label{c1} {Let the conditions of Theorem\/ $\ref{th1}$ be satisfied with  condition~\eqref{5b} replaced by condition~\eqref{3b} with arbitrary
$D\ge\cfrac{1}{2\,\|a\|_{\infty}}$.
Then
$$Q\Big(F_a,\cfrac{1}{D}\Big)\ll \cfrac{1}{\|a\|D\sqrt{M(1)}} +
\exp(-c\,\alpha^2 M(1)). $$}\end{corollary}
\medskip

 It is evident that $M(1)\gg 1-Q(G,2)\geq 1-Q(\mathcal L(X),2)\geq p$, where $p$ is from the conditions of Proposition \ref{thFS}. Note that $M(1)$ can be essentially more than~$p$. 
 For example, $p$ may be equal to $0$, while
$M(1)>0$ for any non-degenerate distribution $F=\mathcal L(X)$.
Therefore, Corollary
\ref{c1} is an essential improvement of Proposition~\ref{thFS}.
It is clear that Corollary \ref{c1} is related to Proposition
\ref{thFS} similarly as  Ess\'een's inequality \eqref{3} is related to the Kolmogorov--Rogozin inequality~\eqref{2}.

 Note that the formulation of Corollary \ref{c1} for each fixed $D$
and for $D=1$ are equivalent. Hence, the formulations of Corollary \ref{c1} 
for all  $D>0$ are equivalent too.

If $D>1$, then ${1}/{D}<1$ and using properties of the concentration function it is easy to see that
Corollary \ref{c1} implies the estimate
 $$Q(F_a,1)\ll D\exp(-c\,\alpha^2 M(1)) + \cfrac{1}{\|a\|\sqrt{M(1)}}\,.$$

The proofs of our Theorem\/ $\ref{th1}$ and Corollary \ref{c1} are in some sence easier than the proofs in 
Friedland and Sodin \cite{Friedland and Sodin} and Rudelson and Vershynin \cite{Rudelson and Vershynin}, since they do not include complicated decompositions of integration sets. This is achieved by an application of relation \eqref{1b}.
Using the methods of Ess\'een \cite{Esseen} (see the proof of Lemma~4 of Chapter~II  in
\cite{Petrov}) is also new in comparison with the arguments in \cite{Friedland and Sodin} 
and~\cite{Rudelson and Vershynin}.
\medskip

Now we reformulate Corollary \ref{c1} for the random variables
${X_k}/{\tau}$, $\tau>0$.

\begin{corollary}\label{c2} {Let
$V_{a,\tau}=\mathcal{L}\big(\sum_{k=1}^{n}a_k {X_k}/{\tau}\big)$. 
Then, under the conditions of Corollary~$\ref{c1}$, we have
\begin{equation}\label{4p}
Q\Big(V_{a,\tau},\cfrac{1}{D}\Big) = Q\Big(F_a,\cfrac{\tau}{D}\Big) \ll
\cfrac{1}{\|a\|D\sqrt{M(\tau)}}
 + \exp(-c\,\alpha^2 M(\tau)).
\end{equation}
Choosing, for example, $\tau=D$, we obtain
$$Q(F_a,1)\ll \cfrac{1}{\|a\|D\sqrt{M(D)}} + \exp(-c\,\alpha^2 M(D)).$$ }
\end{corollary}

 For the proof of Corollary \ref{c2}, it suffices to use  relation \eqref{0}.
 \medskip

 If we consider the special case, where $D=\cfrac{1}{2\,\|a\|_{\infty}}$, then
the restrictions on the arithmetic structure of the vector $a$ are really absent, and we have the bound
\begin{equation}\label{4}
Q(F_a,\tau\,\|a\|_{\infty})\ll \cfrac{\|a\|_{\infty}}{\|a\|\sqrt{M(\tau)}}\,.
\end{equation}
This is just what follows from Ess\'een's inequality applied to the sum of non-identically distributed random variables $Y_k=a_kX_k$ with
$\lambda_k=a_k$,
$\lambda=\|a\|_\infty$.
For $a_1=a_2=\cdots=a_n=1$, inequality \eqref{4} turns into 
the well-known particular case of Proposition \ref{thE}:
\begin{equation}\label{4y}
Q(F^{*n},\tau)\ll \cfrac{1}{\sqrt{n\,M(\tau)}}\,.
\end{equation}
Inequality \eqref{4y} implies also the Kolmogorov--Rogozin inequality for i.i.d. random variables:
$$Q(F^{*n},\tau)\ll \cfrac{1}{\sqrt{n\,(1-Q(F,\tau))}}\,.$$

Inequality \eqref{4} can not give the bound which is better than $O(n^{-1/2})$,
since the right-hand side of \eqref{4} is at least $n^{-1/2}$. 
The results stated above are more interesting if $D$ is essentially more than
~$\cfrac{1}{2\,\|a\|_{\infty}}$. Then one can expect to obtain the estimates 
which are better in order than $O(n^{-1/2})$. Just such estimates of 
$Q(F_a,\lambda)$ are required to study the distributions of eigenvalues of random matrices.

For
$0<D<\cfrac{1}{2\,\|a\|_{\infty}}$ the inequality
\begin{equation}\label{4m}
 Q\Big(F_a,\cfrac{\tau}{D}\Big) \ll
\cfrac{1}{\|a\|D\sqrt{M(\tau)}}
\end{equation}
holds too. In this case it follows from \eqref{8j} and \eqref{4}.

In the statement of Corollary \ref{c2}, the quantity~$\tau$ can be arbitrarily small. If
$\tau$ tends to zero, we obtain
$$Q(F_a,0)\ll \cfrac{1}{\|a\|D\sqrt{{\mathbf P}(\widetilde{X}\neq 0)}} +
\exp(-c\,\alpha^2 \,{\mathbf P}(\widetilde{X}\neq 0)).$$
This estimate could be however deduced  from the results \cite{Friedland and
Sodin} and \cite{Rudelson and Vershynin} too.

\bigskip

Now we formulate improvements of Proposition \ref{thRV} which are similar to Theorem
\ref{th1} and Corollaries \ref{c1} and \ref{c2}.

\begin{theorem}\label{th2}{Let the conditions of Proposition
$\ref{thRV}$ be satisfied for $D=1$.
In other words, $\|\,ta-m\| \geq \min\{\gamma t\,\|a\|,\alpha\}$ 
for all $m=(m_1,\ldots,m_n)\in \mathbf{Z}^n$, $t\in[0,1]$. Then
$$Q(F_a,1) \ll \cfrac{1}{\|a\|\gamma \sqrt{M(1)}} + \exp(-c\,\alpha^2
M(1)).$$}\end{theorem}

\begin{corollary}\label{c3}  {Let the conditions of Proposition
$\ref{thRV}$ be satisfied for an arbitrary~${D>0}$.
Then
$$Q\Big(F_a,\cfrac{1}{D}\Big) \ll \cfrac{1}{\|a\| D \gamma \sqrt{M(1)}} +
\exp(-c\,\alpha^2 M(1)).$$}\end{corollary}

Now we reformulate Corollary \ref{c3} for the vectors ${X_k}/{\tau}$,
$\tau>0$.

\begin{corollary}\label{c4}  {Let
$V_{a,\tau}=\mathcal{L}\big(\sum_{k=1}^{n}a_k {X_k}/{\tau}\big)$.
Then, under the conditions of Corollary~$\ref{c3}$, we have
$$Q\Big(V_{a,\tau},\cfrac{1}{D}\Big) = Q\Big(F_a,\cfrac{\tau}{D}\Big) \ll
\cfrac{1}{\|a\|D\gamma\sqrt{M(\tau)}}
 + \exp(-c\,\alpha^2 M(\tau)).$$
Choosing, for example, $\tau=D$, we have
$$Q(F_a,1)\ll \cfrac{1}{\|a\|D\gamma\sqrt{M(D)}} + \exp(-c\,\alpha^2
M(D)).$$ }
\end{corollary}

 For the proof of Corollary \ref{c4}, it suffices to use  relation \eqref{0}.

\bigskip

 Note that, in Friedland and Sodin \cite{Friedland and Sodin} and Rudelson and
Vershynin \cite{Rudelson and
Vershynin}, the corresponding multi-dimensional results are also contained.
Arguing in a similar way, it is not difficult to transfer the results of this paper to the multivariate case too. In order to simplify the text of this article, we are going to consider the multidimensional case in a different publication.
\bigskip

\section{Proofs}

\emph{Proof of Theorem\/ $\ref{th1}$.} Represent the distribution $G=\mathcal{L}(\widetilde{X})$
as a mixture
$G=q E +\sum_{j=0}^{\infty}p_j G_j $, where $q={\mathbf
P}(\widetilde{X}=0)$, $p_j={\mathbf
P}(\widetilde{X} \in A_j)$, $j=0,1,2,\ldots$, $A_0=\{x:|x|>1\}$,
$A_j=\{x: 2^{-j}<|x|\leq2^{-j+1}\}$, $E$ is probability measure concentrated in zero, $G_j$ 
are probability measures defined for $p_j>0$ by the formula
$G_j\{X\}=\cfrac{1}{p_j}
\,G\{X\bigcap A_j\}$, for any Borel set~$X$. If $p_j=0$,
then we can take as $G_j$ arbitrary measures.

For $z\in \mathbf{R}$, $\gamma>0$, introduce infinitely divisible distributions
$H_{z,\gamma}$, with the characteristic function
$\widehat{H}_{z,\gamma}(t)=\exp\Big(-\cfrac{\gamma}{2}\sum_{k=1}^{n}\big(1-\cos(2a_k
zt)\big)\Big)$.
It is clear that $H_{z,\gamma}$ is a symmetric infinitely divisible distribution. 
It depends on~$a$ too, but we assume that $a$
is fixed.  Therefore, its characteristic function is everywhere positive.

For the characteristic function $\widehat{F}(t)$
of a random variable $X$, we have
$$|\widehat{F}(t)|^2 = \mathbf{E}\exp(it\widetilde{X}) =
\mathbf{E}\cos(t\widetilde{X}),$$
where $\widetilde{X}$ is the corresponding symmetrized random variable. Hence,
\begin{equation}\label{6}|\widehat{F}(t)| \leq
\exp\Big(-\cfrac{\,1\,}{2}\,\big(1-|\widehat{F}(t)|^2\big)\Big)  =
\exp\Big(-\cfrac{\,1\,}{2}\,\mathbf{E}\,\big(1-\cos(t\widetilde{X})\big)\Big).
\end{equation}

 According to \eqref{1} and \eqref{6}, we have
$$Q(F_a,1)\ll \int_{0}^{1}|\widehat{F_a}(t)|\,dt
\ll
\int_{0}^{1}\exp\Big(-\frac{\,1\,}{2}\,\sum_{k=1}^{n}\mathbf{E}\,\big(1-\cos(2a_k
t \widetilde{X})\big)\Big)\,dt=I.$$
It is evident that
\begin{eqnarray*}
\sum_{k=1}^{n}\mathbf{E}\big(1-\cos(2a_k t
\widetilde{X})\big)&=&\sum_{k=1}^{n}\int_{-\infty}^{\infty}\big(1-\cos(2a_k
t x)\big)\,G\{dx\}
 \\
&=&\sum_{k=1}^{n}\sum_{j=0}^{\infty}\int_{-\infty}^{\infty}\big(1-\cos(2a_k
t x)\big)\,p_j
\,G_j\{dx\}\\
&=&\sum_{j=0}^{\infty}\sum_{k=1}^{n}\int_{-\infty}^{\infty}\big(1-\cos(2
a_k t x)\big)\,p_j \,G_j\{dx\}.
\end{eqnarray*}

 We denote $\beta_j=2^{-2j}p_j $,
$\beta=\sum_{j=0}^{\infty}\beta_j$, $\mu_j={\beta_j}/{\beta}$,
$j=0,1,2,\ldots$. It is clear that $\sum_{j=0}^{\infty}\mu_j=1$ and
${p_j}/{\mu_j}=2^{2j}\beta $ (for $p_j> 0$).

 Now we proceed similarly to the proof of a result of Ess\'een \cite{Esseen} (see \cite{Petrov}, Lemma 4 of Chapter II). Using the H\"older inequality, it is easy to see that
 $I\leq \prod _{j=0}^{\infty}I_j^{\mu_j}$, where $I_j=1$ for $p_j=0$. Furthermore, if
$p_j > 0$, then
\begin{eqnarray*}
I_j&=&\int_{0}^{1}\exp\Big(-\cfrac{p_j}{2\,\mu_j}\;\sum_{k=1}^{n}\int_{-\infty}^{\infty}\big(1-\cos(2a_k
t x)\big)\,G_j\{dx\}\Big)\,dt \\
&=&
\int_{0}^{1}\exp\Big(-2^{2j-1}\beta\;\sum_{k=1}^{n}\int_{A_j}\big(1-\cos(2a_k
t
x)\big)\,G_j\{dx\}\Big)\,dt.
\end{eqnarray*}

Applying the Jensen inequality to the exponential in the integral (see
\cite{Petrov},
p. 49)), we obtain
\begin{eqnarray*}
I_j&\leq&\int_{0}^{1}\int_{A_j}\exp\Big(-2^{2j-1}\beta\;
\sum_{k=1}^{n}\big(1-\cos(2a_k t x)\big)\Big)\,G_j\{dx\}\,dt \\
&=
&\int_{A_j}\int_{0}^{1}\exp\Big(-2^{2j-1}\beta\;\sum_{k=1}^{n}\big(1-\cos(2a_k
t x)\big)\Big)\,dt\,G_j\{dx\} \\
&\leq& \sup_{z\in A_j}\int_{0}^{1}\widehat{H}_{z,1}^{2^{2j}\beta}(t)\,dt.
\end{eqnarray*}

Let us estimate the characterictic function $\widehat{H}_{\pi,1}(t)$ for
$|\,t|\leq1$. It is evident that there exists a positive absolute constant $c$ 
such that $1-\cos x \geq cx^2$,
for~${|x|\leq\pi}$. Thus, for $|\,t|\leq \cfrac{1}{2\,\|a\|_{\infty}}$,
\begin{equation}\label{7a}
\widehat{H}_{\pi,1}(t)\leq\exp(-c\|a\|^2t^2).
\end{equation}
For $
\cfrac{1}{2\,\|a\|_{\infty}}\leq|\,t|\leq 1$, one can proceed in the same way as the authors
of~\cite{Friedland and Sodin} and \cite{Rudelson and Vershynin}. Taking into account that $1-\cos t\geq c
\min_{m\in \mathbf{Z}}|\,t-2\pi m|^2$, we obtain
\begin{eqnarray}
\widehat{H}_{\pi,1}(t)&\leq&\exp \Big(-c \;\sum_{k=1}^{n}\min_{m_k \in
\mathbf{Z}}\big|2\pi t a_k -2 \pi m_k\big|^2\Big)\nonumber \\
&=&\exp\Big(-c\;\sum_{k=1}^{n}\min_{m_k \in
\mathbf{Z}}|\,ta_k-m_k|^2\Big)\leq\exp(-c\,\alpha^2),\label{7b}
\end{eqnarray}
for $|\,t|\in \Big[\cfrac{1}{2\,\|a\|_\infty},1\Big]$.

Now we can use inequalities \eqref{7a} and \eqref{7b} to estimate the integrals~$I_j$.
First we consider the case $j=1,2,\ldots$. Note that 
the characteristic functions~$\widehat{H}_{z,\gamma}(t)$ satisfy the equalities
\begin{equation} \label{5}
\widehat{H}_{z,\gamma}(t)=\widehat{H}_{y,\gamma}\big({zt}/{y}\big)\quad\hbox{and}\quad
 \widehat{H}_{z,\gamma}(t)=\widehat{H}_{z,1}^{\gamma}(t).
\end{equation}

For $z\in A_j$ we have $2^{-j}<|z|\leq2^{-j+1}<\pi$. Hence, for
${|\,t|\leq1}$, we have $|{zt}/{\pi}|<1$. Therefore, using the properties
\eqref{5} with $y=\pi$ and aforementioned estimates \eqref{7a} and~\eqref{7b}, we obtain, for $z\in
A_j$,
$$\widehat{H}_{z,1}(t)
\leq\max\big\{\exp\big(-c\,\big({zt\|a\|}/{\pi}\big)^2\big),\;\exp(-c\,\alpha^2)\big\},$$
and, hence,
\begin{eqnarray*}
\sup_{z\in A_j}\int_{0}^{1}\widehat{H}_{z,1}^{2^{2j}\beta}(t)\,dt&\leq&
\int_{0}^{1}\exp(-c\,t^2\beta\|a\|^2)\,dt +
\int_{0}^{1}\exp(-2^{2j}c\,\alpha^2 \beta )\,dt  \\
&\ll &\cfrac{1}{\sqrt{\beta}\,\|a\|} + \exp(-c\,\alpha^2\beta).
\end{eqnarray*}

Consider now the case $j=0$.
The properties \eqref{5} yield, for $z>0,\,\gamma>0$,
\begin{equation}\label{8c}
Q(H_{z,\gamma},1)=Q\big(H_{1,\gamma},{1}/{z}\big).
\end{equation}
Thus, according to  \eqref{8a}, \eqref{1b}, \eqref{5} and
\eqref{8c}, we obtain
\begin{eqnarray*}
\sup_{z\in A_0}\int_{0}^{1}\widehat{H}_{z,1}^{\beta} (t) \,dt &=&
\sup_{z\geq 1} \int_{0}^{1}\widehat{H}_{z,\beta} (t) \,dt \asymp
\sup_{z\geq 1}\; Q(H_{z,\beta},1)\\ &=&
\sup_{z\geq 1}\; Q\big(H_{1,\beta},{1}/{z}\big)
\leq Q(H_{1,\beta},1) \ll Q\big(H_{1,\beta},{1}/{\pi}\big)\\ &=&
Q(H_{\pi,\beta},1) \asymp
\int_{0}^{1} \widehat{H}_{\pi,\beta}(t)\, dt =
\int_{0}^{1}\widehat{H}_{\pi,1}^{\beta}(t) \,dt.
\end{eqnarray*}

Using the bounds \eqref{7a} and \eqref{7b} for the characteristic 
function~$\widehat{H}_{\pi,1}(t)$, we have:
\begin{eqnarray*}
\int_{0}^{1}\widehat{H}_{\pi,1}^\beta(t)\,dt &\leq&
\int_{0}^{1}\exp(-c\|a\|^2\beta t^2)\, dt + \int_{0}^{1}\exp(-c\,\alpha^2
\beta)\,dt \\ &\ll& \cfrac{1}{\|a\|\sqrt{\beta}} + \exp(-c\,\alpha^2 \beta).
\end{eqnarray*}

We obtained the same estimate for all integrals $I_j$ for $p_j\neq 0$. In view of
$\sum_{j=0}^{\infty}\mu_j=1$, we derive that
$$I\leq\prod_{j=0}^{\infty}I_j^{\mu_j} \ll \cfrac{\,1\,}{\|a\|\sqrt{\beta}} +
\exp(-c\,\alpha^2 \beta).$$

Now we will estimate the quantity $\beta$
\begin{eqnarray*}
\beta = \sum_{j=0}^{\infty}\beta_j &=&\sum_{j=0}^{\infty}2^{-2j} p_j  \,
= {\mathbf P}\big(|\widetilde{X}|>1\big) +
\sum_{j=1}^{\infty}2^{-2j}\,{\mathbf
P}\big(2^{-j}<|\widetilde{X}|\leq2^{-j+1}\big)  \\
&\geq&\int_{|x|>1}\,G\{dx\} + \sum_{j=1}^{\infty}
\int_{2^{-j}<|x|\leq2^{-j+1}}\cfrac{x^2}{4}\,G\{dx\}\\ &\geq&
\cfrac{\,1\,}{4}
\int_{|x|>1}\,G\{dx\} + \cfrac{\,1\,}{4} \int_{|x|\leq1}x^2 \,G\{dx\} =
\cfrac{\,1\,}{4} \,M(1).
\end{eqnarray*}
Thus,
\begin{equation}\label{9}
\beta \geq \cfrac{\,1\,}{4}\, M(1).
\end{equation}
Hence,
$$\cfrac{1}{\|a\|\sqrt{\beta}} + \exp(-c\,\alpha^2 \beta) \ll
\cfrac{1}{\|a\|\sqrt{M(1)}} + \exp(-c\,\alpha^2 M(1)),$$
that was required to prove. $\square$

\medskip

Now we will deduce Corollary \ref{c1} from Theorem \ref{th1}.

\emph{Proof of Corollary $\ref{c1}$.} We denote $b=Da \in
\mathbf{R}^n$.
Then the equality $Q(F_a,{1}/{D})=Q(F_b,1)$ is valid. The vector~$b$ satisfies 
the conditions of Theorem \ref{th1} which were there supposed for the vector~$a$. Indeed,
$\|ub-m\|\geq\alpha $ for $u \in \Big[\cfrac{1}{2\,\|b\|_{\infty}},1\Big]$.
This follows from condition~\eqref{3b} of Corollary \ref{c1}, if we denote
$u={t}/{D}$.
It remains to apply Theorem \ref{th1} to the vector $b$. $\square$
\bigskip

\emph{Proof of Theorem\/ $\ref{th2}$.}
We will argue similarly to the proof of Theorem \ref{th1}. Using the notation of Theorem
\ref{th1}, we recall that
$$Q(F_a,1) \ll \prod_{j=0}^{\infty}\sup_{z \in
A_j}\int_{0}^{1}\widehat{H}_{z,1}^{2^{2j}\beta }(t)\,dt \leq
\prod_{j=0}^{\infty}\sup_{z \in
A_j}\int_{0}^{1}\widehat{H}_{\pi,1}^{2^{2j}\beta
}\big({xt}/{\pi}\big)\,dt.$$

The conditions of Theorem \ref{th2} imply that
\begin{eqnarray*}
\widehat{H}_{\pi,1}(t) &\leq
&\exp\Big(-c\;\sum_{k=1}^{n}\min_{m_k\in\mathbf{Z}}\;
\bigl|2\pi t a_k - 2\pi m_k\bigr|^2\Big)\\ &\leq& \exp(-c\,\alpha^2) +
\exp(-c\,t^2\gamma^2\,\|a\|^2)
\end{eqnarray*}
for all  $t\in[0,1]$.
Hence,
\begin{eqnarray*}
Q(F_a,1) &\ll& \int_{0}^{1}\exp(-c\,t^2\gamma^2\beta\, \|a\|^2)\,dt +
\int_{0}^{1}\exp(-c\,\alpha^2\beta)\,dt \\ &\ll&
\cfrac{1}{\gamma\sqrt{\beta}\,\|a\|} + \exp(-c\,\alpha^2 \beta).
\end{eqnarray*}
Now we can use the estimate \eqref{9} for the quantity $\beta$ from the proof of Theorem~\ref{th1}. 
According to this bound, $\beta \geq M(1)/4$.
Then
$$Q(F_a,1) \ll \cfrac{1}{\|a\|\gamma \sqrt{M(1)}} + \exp(-c\,\alpha^2M(1)),$$
that was required to prove. $\square$\medskip

\emph{Proof of Corollary\/ $\ref{c3}$.} This proof is similar to the proof of Corollary~\ref{c1}.
We denote $b=Da \in \mathbf{R}^n$ and $u={t}/{D}$. Then
$\|ub-m\|=\|\,ta-m\|\geq\min\{\gamma t\|a\|,\alpha\}$, for all
$m\in\mathbf{Z}^n$ and $t\in[0,1]$. Thus, the conditions of Theorem~\ref{th2} for the vector~$a$ are valid for the vector~$b$ too. It remains to note that $Q(F_a,{1}/{D})=Q(F_b,1)$ and to apply Theorem  \ref{th2} to the vector $b$.
$\square$

\end{document}